\documentclass[11pt]{article}
\usepackage[english]{babel}
\usepackage{times}      
\usepackage{graphicx}
\usepackage{pstricks}
\usepackage{amsmath,amstext,amssymb,amsfonts}
\allowdisplaybreaks
\newtheorem{lemma}{Lemma}
\newtheorem{theorem}{Theorem}
\newtheorem{corollary}{Corollary}
\newtheorem{definition}{Definition}

\newtheorem{example}{Example}

\newtheorem{remark}{Remark}
\newtheorem{Assumption}{Assumption}

\newcommand{\slim} {\mathop{\rm lim\,sup\,}}
\newcommand{\ilim} {\mathop{\rm lim\,inf\,}}

\def\S{\mathbb{S}}

\def\Y{\mathbb{Y}}
\def\h{\mathbf{I}}
\def\X{\mathbb{X}}

\def\A{\mathbb{A}}

\def\P{\mathbb{P}}
\def\F{\mathbb{F}}

\def\C{\mathbb{C}}

\def\lll{\mathbb{L}}

\def\B{\mathcal{B}}
\def\oo{\mathcal{O}}

\def\fff{\mathtt{F}}

\def\Psii{{\Xi_{\int}}}

\begin{document}

\title{Semi-Uniform
Feller Stochastic Kernels
}

\maketitle

\begin{center}
Eugene~A.~Feinberg \footnote{Department of Applied Mathematics and
Statistics,
 Stony Brook University,
Stony Brook, NY 11794-3600, USA, eugene.feinberg@sunysb.edu},\
Pavlo~O.~Kasyanov\footnote{Institute for Applied System Analysis,
National Technical University of Ukraine ``Igor Sikorsky Kyiv Polytechnic
Institute'', Peremogy ave., 37, build, 35, 03056, Kyiv, Ukraine,\
kasyanov@i.ua.},\ and Michael~Z.~Zgurovsky\footnote{National Technical University of Ukraine
``Igor Sikorsky Kyiv Polytechnic Institute'', Peremogy ave., 37, build, 1, 03056,
Kyiv, Ukraine,\
mzz@kpi.ua
}\\

\bigskip
\end{center}

\begin{abstract}
This paper studies transition probabilities from a Borel subset of a Polish space to a product of two Borel subsets of Polish spaces.  For such transition probabilities it introduces and studies {\color{black}the property of} semi-uniform Feller continuity.  This paper provides several equivalent definitions of semi-uniform Feller continuity {\color{black} and establishes its preservation  under integration}.  The motivation for this study came from the theory of Markov decision processes with incomplete information, and this paper provides fundamental results useful for this theory.
\\
\textbf{Keywords}{ Stochastic kernel; semi-uniform Feller; weak convergence; convergence in total variation; Borel set.}\\
{\bf AMS(2020) subject classification:} Primary 60B10,  Secondary 60J05\\
{\bf Data availability:} Data sharing not applicable to this article as no datasets were generated or analysed during the current study.
\end{abstract}

\section{Introduction}
\label{intro}

This paper studies continuity properties of stochastic kernels, also called transition probabilities, from a Borel subset of a Polish space to a product of two Borel subsets of Polish spaces.  The main property we introduce is semi-uniform Feller continuity, which is a weaker property than continuity in total variation, sometimes called uniform Feller continuity.  This  paper provides several equivalent definitions of semi-uniform Feller continuity.  It also describes the preservation property of semi-uniform Fellerness under integration.

Our main motivation for studying stochastic kernels from a measurable space $\S_3$ to a measurable space $\S_1\times \S_2,$  where $\S_1,$ $\S_2,$ and $\S_3$ are Borel subsets of Polish spaces, is the use of such kernels in mathematical models of decision making with incomplete information. For a Markov decision process with incomplete information, $\S_1$ is an unobservable (or hidden) state space,  $\S_2$ is the set of observations, and $\S_3$ can be either a product of these two spaces and the space of decisions or a subset of the product of these three spaces. Such problems can be reduced to problems with completely observable states by replacing the state space $\S_1$ with the space $\P(\S_1)$ of probability  measures on $\S_1,$ and the new states are called either posterior probabilities or belief states.  This reduction was introduced in \cite{Ao,As,Dy,Shi64,Shi}, and it holds under general {\color{black} measurability} assumptions \cite{Rh,Yu}.

However this reduction does not say much about continuity properties of the transition probability for a new model with the {\color{black} belief} state space $\P(\S_1).$ Weak continuity of this transition probability is essentially necessary for the existence of optimal policies, validity of optimality equations, and convergence of value iterations for models with incomplete information~\cite[Theorem 3.1]{FKZ}.  For models with finite state, observation, and action sets, the required weak continuity of the transition probability in the  model {\color{black} with belief states} takes place~\cite{SS}.

However, weak continuity of the original transition and observation probabilities and even some their stronger properties do not imply weak continuity of the transition probability in the model {\color{black} with belief states}  \cite[Examples 4.1-4.3]{FKZ}. For a partially observable Markov decision process {\color{black} called a ${\rm POMDP}_2$ in this paper and in \cite{FKZSIAM}}, which is a popular particular model of Markov decision process with incomplete information, it was shown in \cite[Theorem 3.6]{FKZ} that weak continuity of transition probabilities and continuity in total of variation of observation probabilities imply  weak continuity of transition probabilities
 in the model with {\color{black} completely observable belief states}
 to which the original problem is reduced.  Another proof of this fact is provided in \cite{Saldi}, where it is also shown that, if the observation probabilities do not depend on controls, then continuity of transition probabilities in total variation imply weak continuity of the transition probabilities in the reduced model with {\color{black} completely observable belief states}.

The remarkable feature of semi-uniform Feller transition kernels is that this property holds in the reduced model with complete information if and only if it holds for the original model, and this fact implies several new and known results on weak continuity of transition probabilities in the {\color{black} model with completely observable states} including all the results described above; see  Section~\ref{secAppl} below and \cite{FKZSIAM} for details.  
This paper provides fundamental results useful for the analysis and optimization of Markov decision processes with complete and incomplete information.  They are used in \cite{FKZSIAM} 
for studying Markov decision processes with incomplete information.  Markov decision processes with semi-uniform Feller transition probabilities are studied in \cite{FKZSIAM} for problems with expected total costs and in \cite{FKZConf}  for problems with average costs per unit time.

For a metric space $\S=(\S,\rho_\S),$ where $\rho_\S$ is a metric, let $\tau(\S)$ be the topology of $\S$ (the family of all open subsets of $\S$), and let ${\mathcal B}(\S)$ be its Borel
$\sigma$-field, that is, the $\sigma$-field generated by all open subsets of the
metric space $\S$. For $s\in \S$ and $\delta>0,$ we denote by
{\color{black}$B(s;\delta):=\{u\in\S:\rho(s,u)< \delta\}$ and $\bar{B}(s;\delta)=\{u\in\S:\rho(s,u)\le \delta\}$ respectively the open and
closed balls in the
metric space $\S$ of  radius $\delta$ with the center $s$ and by $S(s;\delta):=\{u\in\S:\rho(s,u)=\delta\}$ the sphere in
$\S$ of  radius $\delta$ with  center $s$.} 
For a subset $S$ of $\S$ let $\bar{S}$ denote the \textit{closure of} $S,$ and $S^o$ is the \textit{interior of} $S.$ Then $S^o$ is open, $\bar{S}$ is closed, and $S^o\subset S\subset \bar{S}.$ Let $\partial S:=\bar{S}\setminus S^o$ denote the \textit{boundary of} $S.$ {\color{black}We remark that $\partial B(s;\delta)\subset S(s;\delta).$}

We denote by $\P(\S)$ the \textit{set of probability
measures} on $(\S,{\mathcal B}(\S)).$ A sequence of probability
measures $\{\mu^{(n)}\}_{n=1,2,\ldots}$ from $\P(\S)$
\textit{converges weakly} to $\mu\in\P(\S)$ if for any
bounded continuous function $f$ on $\S$
\[\int_\S f(s)\mu^{(n)}(ds)\to \int_\S f(s)\mu(ds) \qquad {\rm as \quad
}n\to\infty.
\]
{\color{black} If this convergence of integrals holds for every bounded  Borel function $f,$  then the sequence $\{\mu^{(n)}\}_{n=1,2,\ldots}$ converges to $\mu$ \emph{setwise.}}
A sequence of probability measures $\{\mu^{(n)}\}_{n=1,2,\ldots}$ from $\P(\S)$ \textit{converges in  total  variation} to $\mu\in\P(\S)$ if
\begin{equation}\label{eq:Kara1}
\begin{aligned}
\sup_{C\in \B(\S)}|\mu^{(n)}(C)-\mu(C)|\to 0\ {\rm as} \ n\to \infty;
\end{aligned}
\end{equation}
see~\cite{Steklov,UFL,Ma} for properties of these types of convergence of probability measures.

Note that $\P(\S)$ is a separable metric space with respect to the topology of weak convergence for probability measures when $\S$ is
a separable metric space; \cite[Chapter~II]{Part}. Moreover,
according to Bogachev \cite[Theorem~8.3.2]{bogachev}, if the metric space $\S$ is separable, then the topology of weak convergence of probability measures on $(\S,\B(\S))$ coincides with the topology generated by the \textit{Kantorovich-Rubinshtein metric}
\begin{equation}\label{eq:KantorRubMetr}
\begin{aligned}
&\rho_{\P(\S)}(\mu,\nu)\\
&:=\sup\left\{\int_\S f(s)\mu(ds)-\int_\S f(s)\nu(ds) \ \Big{|} \ f\in{\rm Lip}_1(\S),\ \sup_{s\in\S}|f(s)|\le 1 \right\},
\end{aligned}
\end{equation}
$\mu,\nu\in\P(\S),$ where
\[
{\rm Lip}_1(\S):=\{f:\S\mapsto\mathbb{R}, \ |f(s_1)-f(s_2)|\le \rho_\S(s_1,s_2),\ \forall s_1,s_2\in\S\}.
\]

For a Borel subset $S$ of a metric space $(\S,\rho_\S)$, where $\rho_\S$
is a metric, we always consider the
 metric space $(S,\rho_S),$ where $\rho_S:=\rho_\S\big|_{S\times S}.$  A subset $B$ of $S$ is called open
(closed) in $S$ if $B$ is open (closed) in $(S,\rho_{\color{black}S})$. Of course, if $S=\S$, we omit
``in $\S$''. Observe that, in general, an open (closed) set
in $S$ may not be open (closed). For $S\in\B(\S)$ we denote by $\B(S)$ the Borel $\sigma$-field on
$(S,\rho_S).$  Observe that $\B(S)=\{S\cap B:B\in\B(\S)\}.$
For metric spaces $\S_1$ and $\S_2$, a (Borel-measurable) \textit{stochastic
kernel} $\Psi(ds_1|s_2)$ on $\S_1$ given $\S_2$ is a mapping $\Psi(\,\cdot\,|\,\cdot\,):\B(\S_1)\times \S_2\mapsto [0,1]$ such that $\Psi(\,\cdot\,|s_2)$ is a
probability measure on $\S_1$ for any $s_2\in \S_2$, and $\Psi(B|\,\cdot\,)$ is a Borel-measurable function on $\S_2$ for any Borel set $B\in\B(\S_1)$.  Another name for a stochastic kernel is a transition probability. A
stochastic kernel $\Psi(ds_1|s_2)$ on $\S_1$ given $\S_2$ defines a Borel measurable mapping $s_2\mapsto \Psi(\,\cdot\,|s_2)$ of $\S_2$ to the metric space
$\P(\S_1)$ endowed with the topology of weak convergence.
A stochastic kernel
$\Psi(ds_1|s_2)$ on $\S_1$ given $\S_2$ is called
\textit{weakly continuous ({\color{black} setwise continuous,} continuous in
total variation)}, if $\Psi(\,\cdot\,|s^{(n)})$ converges weakly ({\color{black}setwise,} in
total  variation) to $\Psi(\,\cdot\,|s)$ whenever $s^{(n)}$ converges to $s$
in $\S_2$. For a singleton $\{s_1\}\subset \S_1,$ we
sometimes write $\Psi(s_1|s_2)$ instead of $\Psi(\{s_1\}|s_2)$. Sometimes a weakly continuous stochastic kernel is called Feller, and a stochastic kernel continuous in total variation is called uniformly Feller \cite{Papa}.

Let $\S_1,$ $\S_2,$ and $\S_3$ be Borel subsets of Polish spaces (a Polish space is a complete
separable metric space), and $\Psi$ on $\S_1\times\S_2$ given $\S_3$ be a stochastic kernel. For $A\in\B(\S_1),$
$B\in\B(\S_2),$ and $s_3\in\S_3,$ let
\begin{equation}\label{eq:marg_new}
\Psi(A,B|s_3):=\Psi(A\times B|s_3).
\end{equation}
In particular, we consider \textit{marginal} stochastic kernels
$\Psi(\S_1,\,\cdot\,|\,\cdot\,)$ on $\S_2$ given $\S_3$ and $\Psi(\,\cdot\,,\S_2|\,\cdot\,)$ on $\S_1$ given $\S_3.$
%
\begin{definition}\label{defi:unifFP}
A stochastic kernel $\Psi$ on $\S_1\times\S_2$ given $\S_3$ is \emph{semi-uniform Feller} if, for each sequence $\{s_3^{(n)}\}_{n=1,2,\ldots}\subset\S_3$ that converges to $s_3$ in $\S_3$ and for each bounded continuous function $f$ on $\S_1,$
\begin{equation}\label{eq:equivWTV3}
\lim_{n\to\infty} \sup_{B\in \B(\S_2)} \left| \int_{\S_1} f(s_1) \Psi(ds_1,B|s_3^{(n)})-\int_{\S_1} f(s_1) \Psi(ds_1,B|s_3)\right|= 0.
\end{equation}
\end{definition}

{\color{black}Definition~\ref{defi:unifFP} implies that for each sequence $\{s_3^{(n)}\}_{n=1,2,\ldots}\subset\S_3$ that converges to $s_3$ in $\S_3,$ for each bounded continuous function $f$ on $\S_1,$ and for each $B\in \B(\S_2),$
\[\lim_{n\to\infty}
 \int_{\S_1} f(s_1) \Psi(ds_1,B|s_3^{(n)})=\int_{\S_1} f(s_1) \Psi(ds_1,B|s_3),
\]
and, in view of Sch\"al \cite[Theorem~3.7(iii,viii)]{sha}, this property implies weak continuity of $\Psi$ on $\S_1\times \S_2$ given $\S_3.$
Thus, a semi-uniform Feller stochastic kernel $\Psi$ on $\S_1\times \S_2$ given $\S_3$ is weakly continuous.}

{\color{black}
The classic definition of weak continuity via convergence of integrals of bounded continuous functions is also applicable to finite measures. It is obvious that a sequence of measures $\{\mu^{(n)}\}_{n=1,2,\ldots}$ on a metric space $\S$ converges weakly to a finite measure $\mu$ on $\S$ if and only if $\mu^{(n)}(\S)\to\mu(\S)$ and, if ${\color{black} \mu}(\S)>0,$ then  the sequence of probability measures  $\{\mu^{(n)}(ds)/\mu^{(n)}(\S)\}_{n=1,2,\ldots}$ converges weakly to the probability {\color{black} measure} $\mu(ds)/\mu(\S).$  In the previous sentence we mean that
$\mu^{(n)}(ds)/\mu^{(n)}(\S)$ is an arbitrary probability measure on $\S$ if $\mu^{(n)}(\S)=0.$

We recall that the marginal measure $\Psi(ds_1,B|s_3),$ $s_3\in\S_3,$ is defined in \eqref{eq:marg_new}. As follows from \eqref{eq:equivWTV3},  if $\Psi$ is a semi-uniform Feller stochastic kernel on $\S_1\times\S_2$ given $\S_3,$ then for each $B\in\mathcal{B}(\S_2)$ the kernel $\Psi(ds_1,B|s_3)$ on $\S_1$ given $\S_3$ is weakly continuous, that is, if $s^{(n)}_3\to s_3$ as $n\to\infty,$ where $s^{(n)}_3,s_3\in\S_3$ for $n=1,2,\ldots,$ then sequence of substochastic measures $\{\Psi(ds_1,B|s^{(n)}_3)\}_{n=1}^\infty$ converges weakly to $\Psi(ds_1,B|s_3).$}
The term ``semi-uniform'' is used in Definition~\ref{defi:unifFP} because {\color{black} the convergence in \eqref{eq:equivWTV3} is uniform only with respect to the second coordinate, and the function $f$ does not depend on the second coordinate}.

This paper describes useful properties of semi-uniform Feller kernels. {\color{black} Section~\ref{sec:PrLESC}, whose main results are Theorem~\ref{th:PrLESC} and its Corollary~\ref{cor:th5.1},  examines the preservation of lower semi-equicontinuity by integrals.  Section~\ref{sec:K}
 studies properties of semi-uniform Feller kernelss.}  Theorem~\ref{th:equivWTV} provides several necessary and sufficient conditions for a stochastic kernel $\Psi$  to be semi-uniform Feller.   Theorem~\ref{th:concept} establishes another necessary and sufficient condition for a stochastic kernel to be semi-uniform Feller. This condition is Assumption~\ref{AssKern}, whose stronger version was introduced in \cite[Theorem~4.4]{Steklov} as a sufficient condition for weak continuity of transition probabilities for Markov decision processes with belief states. Theorem~\ref{th:extra} describes the preservation of {\color{black} semi-uniform Fellerness} under the integration operation.  Section~\ref{secAppl} explains the main  motivation for this study.  Section~\ref{app:B} contains proofs of Theorems~\ref{th:PrLESC}, \ref{th:equivWTV}, \ref{th:concept}, and \ref{th:extra}.   

\section{Preservation of Lower Semi-Equicontinuity by Integrals}\label{sec:PrLESC}


{\color{black}This section provides definitions of equicontinuity properties for families of functions used in this paper and  {\color{black} introduces
Theorem~\ref{th:PrLESC}   stating that integration of the elements of a family of lower semi-equicontinuous functions of two variables in one of these variables preserves lower semi-equicontinuity. Theorem~\ref{th:PrLESC} is used in the proof of Theorem~\ref{th:extra}. }
}

Let us consider some basic definitions.
\begin{definition}\label{defi:semi}
Let $\S$ be a metric space.
 A function $f:\S\to \mathbb{R}$ is called
\begin{itemize}
\item[{\rm(i)}] \textit{lower semi-continuous} (l.s.c.) at a point $s\in\S$ if $\mathop{\ilim}\limits_{s'\to s}f(s')\ge f(s);$
\item[{\rm(ii)}] \textit{upper semi-continuous} at $s\in\S$ if $-f$ is lower semi-continuous at $s;$
\item[{\rm(iii)}]
\textit{continuous} at $s\in\S$ if $f$ is both lower and upper semi-continuous at $s;$
\item[{\rm(iv)}] \textit{lower / upper semi-continuous (continuous respectively) (on $\S$)} if $f$ is lower / upper semi-continuous (continuous respectively) at each $s\in\S.$
\end{itemize}
\end{definition}
For a metric space $\S,$ let  $\F(\S),$ $\lll(\S),$ and $\C(\S)$ be the spaces of all
real-valued functions, all real-valued lower semi-continuous functions, and all real-valued continuous functions respectively defined
on the metric space $\S.$ The following definitions are taken from \cite{FKL2}.

\begin{definition}\label{defi:equi}
A set $\mathtt{F}\subset \F(\S)$ of real-valued functions on a metric space $\S$ is called
\begin{itemize}
\item[{\rm(i)}]
\textit{lower semi-equicontinuous at a point} $s\in \S$ if $\ilim_{s'\to s}\inf_{f\in\mathtt{F}} (f(s')-f(s))\ge0;$
\item[{\rm(ii)}]
\textit{upper semi-equicontinuous at a point} $s\in \S$ if the set $\{-f\,:\,f\in\mathtt{F}\}$ is lower semi-equicontinuous at $s\in \S;$
\item[{\rm(iii)}] \textit{equicontinuous at a point $s\in\S$}, if $\mathtt{F}$ is both lower and upper semi-equiconti-nuous at $s\in\S,$ that is, $\mathop{\lim}\limits_{s'\to s} \mathop{\sup}\limits_{f\in\mathtt{F}} |f(s')-f(s)|=0;$
\item[{\rm(iv)}]
\textit{lower / upper semi-equicontinuous (equicontinuous respectively)} (\textit{on $\S$}) if it is lower / upper semi-equicontinuous (equicontinuous respectively) at all $s \in \S;$
\item[{\rm(v)}] \textit{uniformly bounded (on $\S$)}, if there exists a constant $M<+\infty $ such that $ |f(s)|\le M$ for all $s\in\S$ and  for
all $f\in \mathtt{F}.$
\end{itemize}
\end{definition}

Obviously, if a set
$\mathtt{F}\subset \F(\S)$ is lower semi-equicontinuous, then
$\mathtt{F}\subset \lll(\S).$ Moreover, if a set $\mathtt{F}\subset \F(\S)$ is equicontinuous, then $\mathtt{F}\subset \C(\S).$
The following theorem 
is the main result of this section. 
\begin{theorem}\label{th:PrLESC}
Let $\S_1$, $\S_2$, and $\S_3$ be metric spaces, let $\mathcal{A}\subset \lll(\S_1\times\S_2)$ be a set of functions which is lower semi-equicontinuous and uniformly bounded, and let a stochastic kernel $\Psi(ds_2|s_3)$ on $\S_2$
given $\S_3$ be weakly continuous. If $\S_2$ is separable, then the set of
functions
\begin{equation}\label{eq:FIn}
\mathcal{A}^\Psi:=\left\{(s_1,s_3)\mapsto\int_{\S_2}f(s_1,s_2)\Psi(ds_2|s_3)\,:\, f\in\mathcal{A}\right\}
\end{equation}
defined on $\S_1\times\S_3$ is lower semi-equicontinuous and uniformly bounded by the same constant as the set $\mathcal{A}.$
\end{theorem}

{\color{black}The proof of Theorem~\ref{th:PrLESC} is provided in Section~\ref{app:B}.}

Since $\mathcal{A}\subset \lll(\S_1\times\S_2)$ and $\mathcal{A}$ is uniformly bounded in Theorem~\ref{th:PrLESC}, for each
$s_1\in\S_1$ and $f\in \mathcal{A},$ the bounded function $s_2\mapsto f(s_1,s_2)$ is lower semi-continuous.
Therefore, it is Borel-measurable and bounded. Thus, the integrals in formula \eqref{eq:FIn} are well-defined.

\begin{corollary}\label{cor:th5.1}
Let $\S_1$, $\S_2$, and $\S_3$ be metric spaces, let $\mathcal{A}\subset \mathbb{C}(\S_1\times\S_2)$ be a set of functions which is equicontinuous and uniformly bounded, and let a stochastic kernel $\Psi(ds_2|s_3)$ on $\S_2$
given $\S_3$ be weakly continuous. If $\S_2$ is separable, then the set of
functions ${\mathcal A}^\Psi$ on $\S_1\times\S_3$ defined in \eqref{eq:FIn}
is equicontinuous and uniformly bounded by the same constant as the set $\mathcal{A}.$
\end{corollary}
\underline{Proof.}
{\color{black}Corollary~\ref{cor:th5.1} follows from Theorem~\ref{th:PrLESC} applied to the sets of functions $\mathcal{A}$ and $\{-f\,:\, f\in\mathcal{A}\}.$$~$\hfill$\Box$}
\begin{remark}\label{rem:th5.1}
{\rm
Corollary~\ref{cor:th5.1} is a particular case of  \cite[Theorem 5.1]{FKZ}, where {\color{black} under the same assumption a more general conclusion is stated, which is incorrect.  The difference is that in  \cite[Theorem 5.1]{FKZ} the integration in \eqref {eq:FIn} above is taken over an arbitrary open subset $\mathcal{O}$ of $\S_2$ rather than over the set $\S_2.$  However, the proofs in \cite{FKZ} apply \cite[Theorem 5.1]{FKZ} only to the case $\mathcal{O}=\S_2,$ which is stated in Corollary~\ref{cor:th5.1}.}
}
\end{remark}

{\color{black}Theorem~\ref{th:PrLESC} can be viewed as an extension from equicontinuity to lower semi-equicontinuity of $\mathcal{D}$ {\color{black}in} the following necessary and sufficient condition for  weak convergence of probability measures, whose sufficiency part is obvious by considering a singleton  $\mathcal{D}.$}

\begin{theorem}\label{th:weakconvPar}{\rm(Parthasarathy \cite[Theorem~II.6.8]{Part})}
Let $\S$ be a separable metric space and $(\mu^{(n)})_{n=1,2,\ldots}$
be any sequence of probability measures on $\S.$ Then $\{\mu^{(n)})\}_{n=1,2,\ldots}$
converges weakly to $\mu\in\P(\S)$ if and only if
\[
\lim_{n\to\infty}\sup_{f\in\mathcal{D}}\left|\int_\S f(s)\mu^{(n)}ds-\int_\S f(s)\mu(ds)\right|= 0
\]
for every set $\mathcal{D}\subset \C(\S),$ which is equicontinuous and uniformly bounded.
\end{theorem}

\section{Properties of Semi-Uniform Feller Stochastic Kernels}\label{sec:K}

This section {\color{black}studies the properties of semi-uniform Feller kernels.  In particular, Theorem~\ref{th:equivWTV} provides several necessary and sufficient conditions for semi-uniform Fellerness.}    Theorem~\ref{th:concept} establishes another necessary and sufficient condition for a stochastic kernel to be semi-uniform Feller. This condition is Assumption~\ref{AssKern}, whose stronger version was introduced in \cite[Theorem~4.4]{Steklov}. Theorem~\ref{th:extra} describes the preservation of {\color{black} semi-uniform Feller continuity} under the integration operation.

Let $\S_1,$ $\S_2,$ and $\S_3$ be Borel subsets of Polish spaces, and let $\Psi$ on $\S_1\times\S_2$ given $\S_3$ be a stochastic kernel.
For each set $A\in\B(\S_1)$ consider the set of functions
\begin{equation}\label{eq:familyoffunctions}
\fff^\Psi_A=\{  s_3\mapsto \Psi(A\times B |s_3):\, B\in \B(\S_2)\}
\end{equation}
mapping $\S_3$ into $[0,1].$ Consider the following type of continuity for stochastic kernels on $\S_1\times\S_2$ given $\S_3.$
\begin{definition}\label{defi:wtv}
A stochastic kernel $\Psi$ on $\S_1\times\S_2$ given $\S_3$ is called \textit{WTV-continuous}, if for each $\oo \in\tau (\S_1)$
the set of
functions $\fff^\Psi_\oo$ is lower semi-equicontinuous on $\S_3.$
\end{definition}
 Definition~\ref{defi:equi}{\color{black}(i)} directly implies that the stochastic kernel $\Psi$ on $\S_1\times\S_2$ given $\S_3$ is WTV-continuous if and only if for each $\oo \in\tau (\S_1)$
\begin{equation}\label{eq:equivWTV0}
\ilim_{n\to\infty} \inf_{B\in \B(\S_2)\setminus\{\emptyset\}} \left( \Psi(\oo \times B|s_3^{(n)})-\Psi(\oo \times B|s_3)\right)\ge 0,
\end{equation}
whenever $s_3^{(n)}$ converges to $s_3$ in $\S_3.$  {\color{black}``WTV-continuity'' in Definition~\ref{defi:wtv} abbreviates weak continuity of $\Psi$ in the first variable $s_1\in\S_1$ and continuity in total variation of $\Psi$ in the second variable $s_2\in\S_2.$}

{\color{black} Since} $\emptyset\in\B(\S_2),$ \eqref{eq:equivWTV0} holds if and only if
\begin{equation}\label{eq:wtvsc}
\lim_{n\to\infty} \inf_{B\in \B(\S_2)} \left( \Psi(\oo \times B|s_3^{(n)})-\Psi(\oo \times B|s_3)\right)= 0.
\end{equation}

Similarly to Parthasarathy \cite[Theorem~II.6.1]{Part} {\color{black} and Sch\"al \cite[Theorem~3.7]{sha}}, where necessary and sufficient conditions for weakly convergent probability
measures were considered, the following theorem
provides several useful equivalent definitions of {\color{black}semi-uniform Feller} stochastic kernels.

\begin{theorem}\label{th:equivWTV}
For a stochastic kernel $\Psi$ on $\S_1\times\S_2$ given $\S_3,$  the following conditions are equivalent:
\begin{itemize}
\item[{\rm(a)}] the stochastic kernel $\Psi$ on $\S_1\times\S_2$ given $\S_3$ is semi-uniform Feller;
\item[{\rm(b)}] the stochastic kernel $\Psi$ on $\S_1\times\S_2$ given $\S_3$ is WTV-continuous;
\item[{\rm(c)}] if $s_3^{(n)}$ converges to $s_3$ in $\S_3,$ then for each closed set $C$ in $\S_1$
\begin{equation}\label{eq:wtvscConv}
\lim_{n\to\infty} \sup_{B\in \B(\S_2)} \left( \Psi(C \times B|s_3^{(n)})-\Psi(C \times B|s_3)\right)= 0;
\end{equation}
\item[{\rm(d)}] if $s_3^{(n)}$ converges to $s_3$ in $\S_3,$ then, for each $A\in\B(\S_1)$ such that $\Psi(\partial A,\S_2|s_3)=0,$
\begin{equation}\label{eq:equivWTV2}
\lim_{n\to\infty} \sup_{B\in \B(\S_2)} | \Psi(A \times B|s_3^{(n)})-\Psi(A \times B|s_3)|= 0;
\end{equation}
\item[{\rm(e)}] if $s_3^{(n)}$ converges to $s_3$ in $\S_3,$ then, for each nonnegative bounded lower semi-continuous function $f$ on $\S_1,$
\begin{equation}\label{eq:equivWTV4}
\ilim_{n\to\infty} \inf_{B\in \B(\S_2)} \left( \int_{\S_1} f(s_1) \Psi(ds_1,B|s_3^{(n)})-\int_{\S_1} f(s_1) \Psi(ds_1,B|s_3)\right)= 0;
\end{equation}
\end{itemize}
{\color{black} and each of these conditions implies continuity in total variation of the marginal kernel $\Psi(\S_1,\,\cdot\,|\,\cdot\,)$ on $\S_2$ given $\S_3.$ }
\end{theorem}

The proof of Theorem~\ref{th:equivWTV} is provided in Section~\ref{app:B}.

Note that, since $\emptyset\in\B(\S_2),$ \eqref{eq:wtvscConv} holds if and only if
\begin{equation}\label{eq:equivWTV1}
\slim_{n\to\infty} \sup_{B\in \B(\S_2)\setminus\{\emptyset\}} \left( \Psi(C \times B|s_3^{(n)})-\Psi(C \times B|s_3)\right)\le 0,
\end{equation}
and similar remarks are applicable to \eqref{eq:equivWTV2} and \eqref{eq:equivWTV4} with the inequality ``$\ge$'' taking place in \eqref{eq:equivWTV4}.
%

Let us consider the following assumption. According to Example~\ref{exa:stronger}, Assumption~\ref{AssKern} is weaker than combined assumptions (i) and (ii) in \cite[Theorem~4.4]{Steklov}, where the base $\tau_b^{s_3}(\S_1)$ is the same for all $s_3\in\S_3.$

\begin{Assumption}\label{AssKern}
Let for each $s_3\in\S_3$ the topology on $\S_1$ have a countable base
$\tau_b^{s_3}(\S_1)$ such that
\begin{itemize}
 \item[{\rm(i)}] $\S_1\in\tau_b^{s_3}(\S_1);$
 \item[{\rm(ii)}]  for
each finite intersection $\oo=\cap_{i=1}^ k {\oo}_{i},$ $k=1,2,\ldots,$ of sets
$\oo_{i}\in\tau_b^{s_3}(\S_1),$ $i=1,2,\ldots,k,$
the set of
functions $\fff^\Psi_\oo,$ defined in \eqref{eq:familyoffunctions} with $A=\oo$, is equicontinuous at $s_3.$
\end{itemize}
\end{Assumption}

Note that Assumption~\ref{AssKern}(ii) holds if and only if for
each finite intersection $\oo=\cap_{i=1}^ k {\oo}_{i}$ of sets
$\oo_{i}\in\tau_b^{s_3}(\S_1),$ $i=1,2,\ldots,k,$
\begin{equation}\label{eq:equivWTV0new1}
\lim_{n\to\infty} \sup_{B\in \B(\S_2)} \left| \Psi(\oo \times B|s_3^{(n)})-\Psi(\oo \times B|s_3)\right|= 0
\end{equation}
if $s_3^{(n)}$ converges to $s_3$ in $\S_3.$

The following example demonstrates that the version of Assumption~\ref{AssKern} with the same base $\tau_b(\S_1)$ for all $s_3\in\S_3$ is stronger than Assumption~\ref{AssKern}.

\begin{example}\label{exa:stronger}
{\rm
Let $\S_1=\S_3:=\mathbb{R},$ $\S_2$ be a singleton, and $\Psi(S_1|s_3):=\h\{s_3\in S_1\}$ for all  $S_1\in \B(\S_1)$ and $s_3\in\S_3.$

Let us prove that Assumption~\ref{AssKern} holds. Indeed, for a fixed $s_3\in \mathbb{R}$ let us consider the countable base $\tau_b^{s_3}(\mathbb{R})=\{\mathbb{R}\}\cup\{(a+\sqrt{2},b+\sqrt{2})\,:\,a,b\in \mathbb{Q},\,a<b\}$ for $s_3\in\mathbb{Q},$ and
$\tau_b^{s_3}(\mathbb{R})=\{\mathbb{R}\}\cup\{(a,b)\,:\,a,b\in \mathbb{Q},\,a<b\}$ for $s_3\notin\mathbb{Q},$ where $\mathbb{Q}$ is the set of rational numbers. Note that this base satisfies the following properties:
(a) $\mathbb{R}\in \tau_b^{s_3}(\mathbb{R}),$ (b) $\oo=\cap_{i=1}^ k {\oo}_{i}\in \tau_b^{s_3}(\mathbb{R})$ for any $k=1,2,\ldots$ and
$\{\oo_{i}\}_{i=1}^{k}\subset\tau_b^{s_3}(\mathbb{R}),$ and (c) $s_3\notin\partial \oo$ for all $\oo\in \tau_b^{s_3}(\mathbb{R}).$ Statement~(a) implies that Assumption~\ref{AssKern}(i) holds.
Assumption~\ref{AssKern}(ii) holds because, according to (b) each  finite intersection $\oo=\cap_{i=1}^ k {\oo}_{i}$ of sets
$\oo_{i}\in\tau_b^{s_3}(\mathbb{R}),$ $i=1,2,\ldots,k,$ belongs to $\tau_b^{s_3}(\mathbb{R}),$ and according to (c) the function $s\mapsto \h\{s\in \oo\}$ is continuous at $s_3.$ Thus, Assumption~\ref{AssKern} holds.

Assumption~\ref{AssKern} does not hold with the same base $\tau_b(\S_1)$ for all $s_3\in\S_3$ because for any nonempty open set $\oo\in\tau(\S_1)\setminus \{\S_1\}$ there exist $s_3^*\in \partial \oo$ and a sequence $\{s_3^{(n)}\}_{n=1,2,\ldots}\subset \oo$ such that $s_3^{(n)}\to s_3^*$ in $\S_3$ as $n\to\infty,$ and, therefore, $\Psi(\oo|s_3^{(n)})=\h\{s_3^{(n)}\in \oo\}=1\not\to 0=\h\{s_3^*\in \oo\}=\Psi(\oo|s_3^*)$ as $n\to\infty,$ that is, the set of
functions $\fff^\Psi_\oo$ is not equicontinuous at $s_3^*.$ \hfill $\Box$
}
\end{example}

Theorem~\ref{th:concept} shows that Assumptions~\ref{AssKern} is a necessary and sufficient condition for semi-uniform Feller continuity.

\begin{theorem}\label{th:concept}
{\color{black}A} stochastic kernel $\Psi$ on $\S_1\times\S_2$ given $\S_3$ is semi-uniform Feller if and only if it satisfies Assumption~\ref{AssKern}.
\end{theorem}

The proof of Theorem~\ref{th:concept} is provided in Section~\ref{app:B}.

Now let $\S_4$ be a Borel subset of a Polish space, and let $\Xi$ be a stochastic kernel on $\S_1\times\S_2$ given $\S_3\times\S_4.$ Consider the stochastic kernel $\Psii$ on $\S_1\times\S_2$ given $\P(\S_3)\times\S_4$ defined by
\begin{equation}\label{eq:extra1}
\Psii(A\times B|\mu,s_4):=\int_{\S_3}\Xi(A\times B |s_3,s_4)\mu(ds_3),
\end{equation}
$A\in\B(\S_1),\,B\in\B(\S_2),\,\mu\in\P(\S_3),\,s_4\in\S_4.$

{\color{black}Note that $\Xi$ is the integrand for $\Psii,$ which justifies the notation $\Psii.$}
The following theorem establishes the preservation of {\color{black}semi-uniform Fellerness under} the integration  operation in \eqref{eq:extra1}.

\begin{theorem}\label{th:extra}
{\color{black}A} stochastic kernel $\Psii$ on
$\S_1\times\S_2$ given $\P(\S_3)\times\S_4$ is {\color{black}semi-uniform Feller}
if and only if \, $\Xi$ on $\S_1\times\S_2$ given $\S_3\times\S_4$ is {\color{black}semi-uniform Feller}.
\end{theorem}

The proof of Theorem~\ref{th:extra} is provided in Section~\ref{app:B}.

{\color{black}
\section{Motivation for Studying  Semi-Uniform Feller Continuity: Control of Markov Processes with Incomplete Information}\label{secAppl}

{\color{black}Semi-uniform Feller continuity appears naturally in control of stochastic processes with incomplete information, when a decision maker observes random variables depending on the states of the process rather than the states themselves. The main approach to  analyzing such problems is to consider a stochastic process whose states are posterior distributions of the states of the original process; see, e.g., \cite{Rh,Yu}. These posterior distributions are often called beliefs or belief states.  The Bayesian approach to the control of stochastic processes is based on substituting states of the process with their beliefs.

{\color{black}Let $\X,$ $\Y,$ and $\A$ be Borel subsets of Polish spaces, where $\X$ is the set of hidden states, $\Y$ is the set of observations, and $\A$ is the sets of controls.  Let $P$ be a stochastic kernel on $\X\times\Y$ given $\X\times\Y\times\A.$}
The dynamics of a Markov Decision Process with Incomplete Information (MDPII) \cite{DY,FKZSIAM} is defined {\color{black} by} $P(dx_{t+1},dy_{t+1}|x_t,y_t,a_t),$ where $x_t$ is a hidden state, $y_t$ is an observation, and $a_t$ is a chosen control, $t=0,1,\ldots.$  {\color{black}It} is possible to construct a completely observable Markov Decision Process (MDP) whose dynamics is defined by a stochastic kernel $q(dz_{t+1},dy_{t+1}|z_t,y_t,a_t),$ where $z_t$ is a posterior probability distribution of the state $x_t,$ $t=0,1,\ldots,$ and $q$ can be constructed from $P$ by using the Bayesian arguments  \cite{DY, FKZSIAM, Rh, Yu}.

An important question is whether the transition kernel $q$ is weakly continuous, and weak continuity of kernels is sometimes called Feller continuity. It is known that weak continuity of $P$ does not imply weak continuity of $q$ \cite[Example 4.1]{FKZ}, and finding sufficient conditions for weak continuity of $q$ is an important question.  According to \cite[Theorem 6.2]{FKZSIAM}, whose proof uses the results of this paper, $q$ is semi-uniform Feller if and only if $P$ is semi-uniform Feller.  Thus, semi-uniform Feller continuity of $P$ is a natural sufficient condition for weak continuity of $q.$

Weak continuity of the stochastic kernel $q$ implies weak continuity of its marginal kernel $\hat{q}(dz_{t+1}|z_t,y_t,a_t){\color{black}:=q(dz_{t+1},\Y|z_t,y_t,a_t)}.$ An important particular case of a MDPII is a Partially Observable Markov Decision Process (POMDP). For a POMDP the kernel $P$ has a special structure, which is not important here, but it is important that transition probabilities defined by kernels $P$ do not depend on observations.  This means that the transition probabilities $P$ and $q$ in the case of an POMDP can be written as $P(dx_{t+1},dy_{t+1}|x_t,a_t),$ and $\hat{q}(dz_{t+1}|z_t,a_t).$  If cost functions also do not depend on observations, then the information about observation is useless for the model with belief states constructed for the POMDP.  In this case, the central question is weak continuity of $\hat{q}(dz_{t+1}|z_t,a_t).$

In nonlinear filtering theory, weak continuity of $\hat{q}(dz_{t+1}|z_t,a_t)$ is called weak continuity of the filter \cite{Saldi}. Sufficient conditions for continuity of nonlinear filters and a slightly more general problem of weak continuity of the stochastic kernel $\hat{q}(dz_{t+1}|z_t,a_t)$  for POMDPs were studied recently in \cite{FKZ, FKZSIAM, Saldi}. Earlier results can be found in \cite{HL} and \cite{RS}.   All currently known sufficient conditions for weak continuity of $\hat{q}$ assume semi-uniform Feller continuity of the stochastic kernel $P;$ see \cite[Corollaries 6.10 and 6.11]{FKZSIAM} and \cite{FKZ22}.}

\section{Proofs of Theorems~{\color{black}\ref{th:PrLESC},} \ref{th:equivWTV}, \ref{th:concept}, and \ref{th:extra}}\label{app:B}



Before {\color{black} proving Theorem~\ref{th:PrLESC}} we provide additional definitions and establish additional properties of functions from $\lll(\S_1\times\S_2).$  For a bounded function $g$ defined on a metric space $\S${\color{black},
let us consider its Pasch-Hausdorff envelope defined for $m=1,2,\ldots,$}
\begin{equation}\label{eq:appdefiop1}
r_{g(\,\cdot\,)}^{(m)}(s):=\inf_{s'\in \S} [g(s')+m \rho_{\S}(s,s')],\quad s\in \S;
\end{equation}
{\color{black}see Bertsekas and Shreve \cite[p.~125]{BS}, Rockafellar and Wets \cite{RW}, and Feinberg et al \cite{FKR} for properties of functions defined in \eqref{eq:appdefiop1}.}
Formula \eqref{eq:appdefiop} below defines {\color{black} a parameterized version of the  Pasch-Hausdorff envelope defined for a bounded function $f$ on $\S_1\times \S_2,$} where  the variable $s_1$ plays the role of a parameter, and the variable $s_2$ plays the role of the variable $s$ in \eqref{eq:appdefiop1}. For each $m=1,2,\ldots,$ and $s_1\in\S_1,$ we set
\begin{equation}\label{eq:appdefiop}
r_{f(s_1,\,\cdot\,)}^{(m)}(s_2):=\inf_{s_2'\in \S_2} [f(s_1,s_2')+m \rho_{\S_2}(s_2,s_2')],\quad s_2\in \S_2.
\end{equation}
Let the set of functions $\mathcal{A}$ from Theorem~\ref{th:PrLESC} be uniformly bounded by a constant $M.$ According to Bertsekas and Shreve \cite[p.~125]{BS}, for each $f\in\mathcal{A},$ $m_1,m_2=1,2,\ldots, $ $m_1\le m_2,$ $s_1\in\S_1,$ and $s_2\in \S_2,$ the following inequalities hold,
    \begin{equation}\label{eq:app3}
-M\le    r_{f(s_1,\,\cdot\,)}^{(m_1)}(s_2)\le r_{f(s_1,\,\cdot\,)}^{(m_2)}(s_2)\le f(s_1,s_2).
    \end{equation}

For each $m=1,2,\ldots$ we set
\begin{equation}\label{eq:app2}
\mathcal{C}(\mathcal{A},m):=\{s_2\mapsto r_{f(s_1,\,\cdot\,)}^{(m)}(s_2)\,:\, f\in\mathcal{A},\, s_1\in\S_1\}\subset \F(\S_2).
\end{equation}
The following lemma establishes basic properties of  the sets $\mathcal{C}(\mathcal{A},m),$ $m=1,2,\ldots\,.$ It is used in the proofs of Theorems~\ref{th:PrLESC} and \ref{th:equivWTV}.  It describes uniform approximations of functions in families of lower semi-continuous functions by globally Lipschitz functions.
\begin{lemma}\label{lem:A1}
Let $\mathcal{A}\subset \lll(\S_1\times\S_2),$ where $\S_1$ and $\S_2$ are metric spaces.  The following statements hold:
\begin{itemize}
\item[{\rm(i)}] if the set $\mathcal{A}$ is  uniformly bounded by a constant $M>0,$ then for each $m=1,2,\ldots$  the set $\mathcal{C}(\mathcal{A},m)$ defined in \eqref{eq:app2} is uniformly bounded by the same constant $M;$
\item[{\rm(ii)}] for each $m=1,2,\ldots$ the set $\mathcal{C}(\mathcal{A},m)$ is equicontinuous;
\item[{\rm(iii)}] if $\mathcal{A}$ is lower semi-equicontinuous and uniformly bounded, then, for each sequence $\{s_1^{(n)}\}_{n=1,2,\ldots}\subset \S_1$ that converges to $s_1\in \S_1$ and for each $s_2\in\S_2,$
\begin{equation}\label{eq:app5}
\ilim_{m\to\infty}\ilim_{n\to\infty}\inf_{f\in \mathcal{A}}\,[r_{f(s_1^{(n)},\,\cdot\,)}^{(m)}(s_2)-f(s_1,s_2)]\ge 0.
\end{equation}
\end{itemize}
\end{lemma}

Lemma~\ref{lem:A1}(iii) is relevant to Bertsekas and Shreve \cite[Lemma~7.14(a)]{BS} stating how a lower semi-continuous function can be approximated from below by continuous functions. If $\mathcal{A}$ consists of one function $f\in\lll(\S_2),$ which does not depend on $s_1,$ then Lemma~\ref{lem:A1}(iii) implies that $r_{f}^{(m)}(s_2)\uparrow f(s_2)$ as $m\to\infty$ for each $s_2\in\S_2$ because $r_{f}^{(m_1)}(s_2)\le r_{f}^{(m_2)}(s_2)\le f(s_2),$
for each $s_2\in\S_2$ and for all $m_1,m_2=1,2,\ldots$ such that $m_1\le m_2.$ Therefore, \eqref{eq:app5} transforms to $0\le f(s_2)- r_{f}^{(m)}(s_2) \downarrow 0$ as $m\to\infty,$ which is equivalent to the conclusion of \cite[Lemma~7.14(a)]{BS} stating that $r_{f}^{(m)}(s_2)\uparrow f(s_2)$ as $m\to\infty$ for each $s_2\in\S_2.$

\underline{Proof of Lemma~\ref{lem:A1}.} (i) According to \eqref{eq:app3}, the set $\mathcal{C}(\mathcal{A},m),$ $m=1,2,\ldots,$ is uniformly bounded by $M$ whenever the set $\mathcal{A}$ is  uniformly bounded by $M.$

(ii) According to Bertsekas and Shreve \cite[pp.~125, 126]{BS}, for each $m=1,2,\ldots, $ $f\in\mathcal{A},$ $s_1\in\S_1,$ and $s_2^{(1)},s_2^{(2)}\in \S_2,$
\begin{equation}\label{eq:app4}
|r_{f(s_1,\,\cdot\,)}^{(m)}(s_2^{(1)})-r_{f(s_1,\,\cdot\,)}^{(m)}(s_2^{(2)})|\le m\rho_{\S_2}(s_2^{(1)},s_2^{(2)}).
\end{equation}
Therefore, for each $m=1,2,\ldots$ the set $\mathcal{C}(\mathcal{A},m)$ is equicontinuous.

(iii) Since $\mathcal{A}$ is uniformly bounded by a constant $M>0,$
\begin{equation}\label{eq:app4_1}
\sup_{f\in\mathcal{A}}\sup_{u_1\in\S_1,\,u_2\in\S_2} |f(u_1,u_2)|\le M.
\end{equation}
Let $m=1,2,\ldots\,,$ $s_i\in \S_i$ for $i=1,2,$ and let us fix an arbitrary sequence $\{s_1^{(n)}\}_{n=1,2,\ldots}\subset \S_1$ converging to $s_1.$  Inequalities \eqref{eq:app3} and \eqref{eq:app4_1} imply that
\begin{equation}\label{eq:app6}
-\infty<-2M\le r_{f(s_1^{(n)},\,\cdot\,)}^{(m)}(s_2)-f(s_1,s_2)\le 2 M<\infty
\end{equation}
for each $f\in\mathcal{A}$ and for an arbitrary integer $n\ge m.$ Let us take the infimum in $n\ge m$ and in $f\in\mathcal{A}$ of the central expression in \eqref{eq:app6}. Since the infimum in two parameters is equal to the double infimum, the definition of an infimum implies the existence of an integer $n(m)\ge m$ and a function $f^{(m)}\in\mathcal{A}$ such that
\begin{equation}\label{eq:app7}
\begin{aligned}
\inf_{n=m,m+1,\ldots}\inf_{f\in \mathcal{A}}&\,[r_{f(s_1^{(n)},\,\cdot\,)}^{(m)}(s_2)-f(s_1,s_2)]\\ &> r_{f^{(m)}(s_1^{(n(m))},\,\cdot\,)}^{(m)}(s_2)-f^{(m)}(s_1,s_2) - \frac1m.
\end{aligned}
\end{equation}
Note that
\begin{equation}\label{eq:app6fr12}
s_1^{(n(m))}\to s_1 \quad\mbox{as}\quad m\to\infty.
\end{equation}

Statement~(i) and formula \eqref{eq:app4_1} imply that, for all $g\in \cal{A}$ and $u\in \S_1,$
\begin{equation}\label{eq:app7_1}
|r^{(m)}_{g(u,\,\cdot\,)}(s_2)|\le M.
\end{equation}
Therefore, $r_{f^{(m)}(s_1^{(n(m))},\,\cdot\,)}^{(m)}(s_2)$ is bounded by $M,$ and, in virtue of \eqref{eq:appdefiop}, there exists $s_2^{(m)}\in \S_2$ such that
\begin{equation}\label{eq:app8}
r_{f(s_1^{(n(m))},\,\cdot\,)}^{(m)}(s_2) > f^{(m)}(s_1^{(n(m))},s_2^{(m)})+m \rho_{\S_2}(s_2,s_2^{(m)})-\frac1m.
\end{equation}

Inequalities \eqref{eq:app8}, \eqref{eq:app4_1} and \eqref{eq:app7_1} imply
$\rho_{\S_2}(s_2,s_2^{(m)})\le \frac{2M}{m}+\frac1{m^2}.$ Therefore,
\begin{equation}\label{eq:app91}
s_2^{(m)}\to s_2 \quad \mbox{as}\quad m\to\infty.
\end{equation}
Inequalities \eqref{eq:app7} and \eqref{eq:app8} imply
\begin{equation}\label{eq:app9}
\begin{aligned}
\inf_{n=m,m+1,\ldots}&\inf_{f\in \mathcal{A}} \,[r_{f(s_1^{(n)},\,\cdot\,)}^{(m)}(s_2)-f(s_1,s_2)]\\
&> f^{(m)}(s_1^{(n(m))},s_2^{(m)})-f^{(m)}(s_1,s_2)+m \rho_{\S_2}(s_2,s_2^{(m)})- \frac2m.
\end{aligned}
\end{equation}

Since $m=1,2,\ldots$ is arbitrary,
\[
\begin{aligned}
\ilim_{m\to\infty}&\ilim_{n\to\infty}\inf_{f\in \mathcal{A}}\,[r_{f(s_1^{(n)},\,\cdot\,)}^{(m)}(s_2)-f(s_1,s_2)]\\
&\ge \ilim_{m\to\infty}\inf_{n=m,m+1,\ldots}\inf_{f\in \mathcal{A}}\,[r_{f(s_1^{(n)},\,\cdot\,)}^{(m)}(s_2)-f(s_1,s_2)]\\
&\ge \ilim_{m\to\infty} [f^{(m)}(s_1^{(n(m))},s_2^{(m)})-f^{(m)}(s_1,s_2)]\\ &\ge  \ilim_{m\to\infty}\inf_{g\in\mathcal{A}} [g(s_1^{(n(m))},s_2^{(m)})-g(s_1,s_2)]\ge0,
\end{aligned}
\]
where the first inequality holds because the lower limit of a sequence is greater than or equal to its infimum; the second inequality follows from \eqref{eq:app9};
the third inequality holds because $\{f^{(m)}\}_{m=1,2,\ldots}\subset \mathcal{A};$ and the last inequality holds because the set $\mathcal{A}$ is lower semi-equicontinuous and because of \eqref{eq:app6fr12} and \eqref{eq:app91}.$~$\hfill$\Box$

\underline{Proof of Theorem~\ref{th:PrLESC}.} Since $\Psi(ds_2|s_3)$ is a stochastic kernel, and since the set of functions $\mathcal{A}\subset \lll(\S_1\times\S_2)$ is
uniformly bounded, the set of
functions $\mathcal{A}^\Psi$ is uniformly bounded by the same constant as $\mathcal{A}$.

Let us prove that the set of
functions $\mathcal{A}^\Psi$ is lower semi-equicontinuous. Fix an arbitrary
sequence $\{s_1^{(n)},s_3^{(n)}\}_{n=1,2,\ldots}\subset \S_1\times\S_3,$ that converges to some  $(s_1,s_3)\in \S_1\times\S_3,$ and fix an arbitrary $m=1,2,\ldots\,.$ Let us define
\[
\begin{aligned}
I_1^{(m)}:=& \ilim_{n\to\infty} \inf_{f\in\mathcal{A}}\left(\int_{\S_2}r_{f(s_1^{(n)},\,\cdot\,)}^{(m)}(s_2)\Psi(ds_2|s_3^{(n)})- \int_{\S_2}r_{f(s_1^{(n)},\,\cdot\,)}^{(m)}(s_2)\Psi(ds_2|s_3)\right),\\
I_2^{(m)}:=&\ilim_{n\to\infty} \inf_{f\in\mathcal{A}}\int_{\S_2}[r_{f(s_1^{(n)},\,\cdot\,)}^{(m)}(s_2)- f(s_1,s_2)]\Psi(ds_2|s_3).
\end{aligned}
\]
Then
\begin{equation}\label{eq:PrLESC2}
\begin{aligned}
&\ilim_{n\to\infty} \inf_{f\in\mathcal{A}}\left(\int_{\S_2}f(s_1^{(n)},s_2)\Psi(ds_2|s_3^{(n)})- \int_{\S_2}f(s_1,s_2)\Psi(ds_2|s_3)\right)\\ \ge&\ilim_{n\to\infty} \inf_{f\in\mathcal{A}}\left(\int_{\S_2}r_{f(s_1^{(n)},\,\cdot\,)}^{(m)}(s_2)\Psi(ds_2|s_3^{(n)})- \int_{\S_2}f(s_1,s_2)\Psi(ds_2|s_3)\right)\\ &\ge I_1^{(m)}+I_2^{(m)},
\end{aligned}
\end{equation}
where the first inequality follows from the last inequality in \eqref{eq:app3}, and the second inequality follows from the semiadditive properties of infimums and lower limits.
Theorem~\ref{th:weakconvPar}, applied to
$\S:=\S_2,$ $\mathcal{D}^{(m)}:=\{r_{f(s_1^{(n)},\,\cdot\,)}^{(m)}\,:\, f\in\mathcal{A},\,n=1,2,\ldots\},$ $\mu^{(n)}(ds_2):=\Psi(ds_2|s_3^{(n)}),$ $n=1,2,\ldots,$ and $\mu(ds_2):=\Psi(ds_2|s_3),$ implies
\begin{equation}\label{eq:PrLESC3}
I_1^{(m)}\ge 0
\end{equation}
because, according to Lemma~\ref{lem:A1}(i,ii), the set of functions $\mathcal{D}^{(m)}\subset \C(\S_2)$ is equicontinuous and uniformly bounded.

Since the sets of functions $\mathcal{D}^{(m)}\subset \C(\S_2)$ and $\mathcal{A}\subset \lll(\S_1\times\S_2)$ are uniformly bounded, the function $s_2\mapsto \mathop{\inf}\limits_{f\in\mathcal{A}}[r_{f(s_1^{(n)},\,\cdot\,)}^{(m)}(s_2)- f(s_1,s_2)]$ is bounded, and it is upper semi-continuous as an infimum of upper semi-continuous functions. Thus, this function is Borel-measurable. Therefore,
\[
\begin{aligned}
I_2^{(m)}&\ge   \ilim_{n\to\infty} \int_{\S_2}\inf_{f\in\mathcal{A}}[r_{f(s_1^{(n)},\,\cdot\,)}^{(m)}(s_2)- f(s_1,s_2)]\Psi(ds_2|s_3)\\
&\ge \int_{\S_2}\ilim_{n\to\infty} \inf_{f\in\mathcal{A}}[r_{f(s_1^{(n)},\,\cdot\,)}^{(m)}(s_2)- f(s_1,s_2)]\Psi(ds_2|s_3),
\end{aligned}
\]
where the first inequality is obvious, and
the second one follows from Fatou's lemma because, according to Lemma~\ref{lem:A1}(i), the set of functions\\ $\{s_2\mapsto \mathop{\inf}\limits_{f\in\mathcal{A}}[r_{f(s_1^{(n)},\,\cdot\,)}^{(m)}(s_2)- f(s_1,s_2)]\}_{n,m=1,2,\ldots}$ is uniformly bounded. Furthermore,
\begin{equation}\label{eq:PrLESCstar}
\begin{aligned}
\ilim_{m\to\infty}& I_2^{(m)}\\ &\ge \int_{\S_2}\ilim_{m\to\infty}\ilim_{n\to\infty} \inf_{f\in\mathcal{A}}[r_{f(s_1^{(n)},\,\cdot\,)}^{(m)}(s_2)- f(s_1,s_2)]\Psi(ds_2|s_3)\ge0,
\end{aligned}
\end{equation}
where  the first inequality follows from Fatou's lemma because the functions to which Fatou's lemma is applied are uniformly bounded in view of Lemma~\ref{lem:A1}(i), and the second inequality follows from Lemma~\ref{lem:A1}(iii). Inequalities \eqref{eq:PrLESC2}, \eqref{eq:PrLESC3}, and \eqref{eq:PrLESCstar} imply 
\[
\begin{aligned}
&\ilim_{n\to\infty} \inf_{f\in\mathcal{A}}\left(\int_{\S_2}f(s_1^{(n)},s_2)\Psi(ds_2|s_3^{(n)})- \int_{\S_2}f(s_1,s_2)\Psi(ds_2|s_3)\right)\\
&\ge \liminf_{m\to\infty} (I_1^{(m)}+I_2^{(m)})\ge  \liminf_{m\to\infty} I_2^{(m)}\ge 0,
\end{aligned}
\]
that is, the set of
functions $\mathcal{A}^\Psi$ is lower semi-equicontinuous.$~$\hfill$\Box$

\underline{Proof of Theorem~\ref{th:equivWTV}.}   {\color{black} Under each of conditions (a)--(e) the marginal kernel\\ $\Psi(\S_1,\,\cdot\,|\,\cdot\,)$ on $\S_2$ given $\S_3$ is continuous in total variation. In particular,} {\color{black} under condition (a) this follows from \eqref{eq:equivWTV3} with $f\equiv 1.$  Under condition (b), continuity in total variation of the marginal kernel $\Psi(\S_1,\,\cdot\,|\,\cdot\,)$ follows from
\begin{equation*}
\begin{aligned}
&\lim_{n\to\infty} \sup_{B\in \B(\S_2)} \left| \Psi(\S_1 \times B|s_3^{(n)})-\Psi(\S_1 \times B|s_3)\right|\\&=\lim_{n\to\infty} \sup_{B\in \B(\S_2)} \left( \Psi(\S_1 \times B|s_3^{(n)})-\Psi(\S_1 \times B|s_3)\right)= 0,
\end{aligned}
\end{equation*}
where the second equality
follows from equality \eqref{eq:wtvsc} with $\oo:=\S_1$ and from $\Psi(\S_1\times\S_2|\,\cdot\,)=1.$ Conditions (c) and (d) with $C=\S_1$ and $A=\S_1$ respectively imply continuity in total variation of this marginal kernel.  In addition, condition (e) with $f(s_1)={\bf I}\{s_1\in\oo\},$ where $\oo$ are open subsets of $\S_1,$ implies condition (b).}

{\color{black} The equivalence of conditions (a)--(e) follows the following implications: $\text{(a)}\, \Rightarrow\, \text{(e)}\, \Rightarrow \text{(b)}\, \Leftrightarrow\, \text{(c)}\,  \Rightarrow\, \text{(d)}\,   \Rightarrow\, \text{(a)}.$}

(a) $\Rightarrow$ (e). Let $s_3^{(n)}$ converge to $s_3$ in $\S_3,$ and let $f$ be a nonnegative bounded lower semi-continuous function on $\S_1.$
We shall prove \eqref{eq:equivWTV4}. Indeed, for an arbitrary fixed $m=1,2,\ldots,$ in view of \eqref{eq:appdefiop1} and the last inequality in \eqref{eq:app3},
\begin{equation}\label{eq:equivWTV5}
\begin{aligned}
\ilim_{n\to\infty}\inf_{B\in \B(\S_2)} &\left( \int_{\S_1} f(s_1) \Psi(ds_1,B|s_3^{(n)})-\int_{\S_1} f(s_1) \Psi(ds_1,B|s_3)\right)\\
\ge\ilim_{n\to\infty} \inf_{B\in \B(\S_2)}& \left( \int_{\S_1} r_{f(\cdot)}^{(m)}(s_1) \Psi(ds_1,B|s_3^{(n)}) \right.\\ &\left.-\int_{\S_1} f(s_1) \Psi(ds_1,B|s_3)\right)\ge I_1^{(m)}+I_2^{(m)},
\end{aligned}
\end{equation}
where
\begin{equation}\label{eq:equivWTV6}
\begin{aligned}
I_1^{(m)}:=\ilim_{n\to\infty} \inf_{B\in \B(\S_2)}& \left( \int_{\S_1} r_{f(\cdot)}^{(m)}(s_1) \Psi(ds_1,B|s_3^{(n)}) \right. \\ &\left.-\int_{\S_1} r_{f(\cdot)}^{(m)}(s_1) \Psi(ds_1,B|s_3)\right)= 0,
\end{aligned}
\end{equation}
\begin{equation}\label{eq:equivWTV7}
\begin{aligned}
I_2^{(m)}:=&\inf_{B\in \B(\S_2)} \left( \int_{\S_1} r_{f(\cdot)}^{(m)}(s_1) \Psi(ds_1,B|s_3)-\int_{\S_1} f(s_1) \Psi(ds_1,B|s_3)\right)\\
&=\int_{\S_1} \left(r_{f(\cdot)}^{(m)}(s_1)-f(s_1)\right) \Psi(ds_1,\S_2|s_3).
\end{aligned}
\end{equation}
We note that the last equality in \eqref{eq:equivWTV6} follows from statement~({a}) because, according to Lemma~\ref{lem:A1}(i,ii), the function $s_1\mapsto r_{f(\cdot)}^{(m)}(s_1)$ is continuous and bounded on $\S_1,$ and the last equality in \eqref{eq:equivWTV7} follows from the inequality $r_{f(\cdot)}^{(m)}(s_1)\le f(s_1)$ for each $s_1\in\S_1.$ Finally, \eqref{eq:equivWTV5}--\eqref{eq:equivWTV7} imply that for each $m=1,2,\ldots$
\[
\begin{aligned}
\ilim_{n\to\infty} &\inf_{B\in \B(\S_2)} \left( \int_{\S_1} f(s_1) \Psi(ds_1,B|s_3^{(n)})-\int_{\S_1} f(s_1) \Psi(ds_1,B|s_3)\right)\\
&\ge \int_{\S_1} \left(r_{f(\cdot)}^{(m)}(s_1)-f(s_1)\right) \Psi(ds_1,\S_2|s_3)\to 0,\quad m\to\infty,
\end{aligned}
\]
where the convergence to zero directly follows from Lebesgue's dominated convergence theorem  because, according to Lemma~\ref{lem:A1}, the sequence $\{r_{f(\cdot)}^{(m)}(\,\cdot\,)-f(\,\cdot\,)\}_{m=1,2,\ldots}$ is uniformly bounded and converges pointwise to zero. Thus, \eqref{eq:equivWTV4} holds.

(e) $\Rightarrow$ (b). Let $s_3^{(n)}$ converge to $s_3$ in $\S_3,$ and $\oo \in\tau (\S_1).$
For a nonnegative bounded lower semi-continuous function $f(s_1):=\h\{s_1\in \oo\},$ $s_1\in\S_1,$ \eqref{eq:equivWTV4} directly implies \eqref{eq:wtvsc} and therefore \eqref{eq:equivWTV0}.

(b) $\Leftrightarrow$ (c). Let $s_3^{(n)}$ converges to $s_3$ in $\S_3.$ Note that for each
$S\in\B(\S_1)$
\[
\begin{aligned}
\slim_{n\to\infty}& \sup_{B\in \B(\S_2)\setminus\{\emptyset\}} \left( \Psi(S \times B|s_3^{(n)})-\Psi(S \times B|s_3)\right)\\
=\slim_{n\to\infty}& \sup_{B\in \B(\S_2)\setminus\{\emptyset\}} \left( \Psi(\S_1\times B|s_3^{(n)})-\Psi(\S_1\times B|s_3)\right. \\ & \left.- \Psi((\S_1\setminus S) \times B|s_3^{(n)})+\Psi((\S_1\setminus S) \times B|s_3)\right) \\
\le \slim_{n\to\infty}& \sup_{B\in \B(\S_2)\setminus\{\emptyset\}} \left| \Psi(\S_1,B|s_3^{(n)})-\Psi(\S_1,B|s_3) \right| \\
-\ilim_{n\to\infty}& \inf_{B\in \B(\S_2)\setminus\{\emptyset\}} \left( \Psi((\S_1\setminus S) \times B|s_3^{(n)})-\Psi((\S_1\setminus S) \times B|s_3)\right)\\  = -\ilim_{n\to\infty}& \inf_{B\in \B(\S_2)\setminus\{\emptyset\}} \left( \Psi((\S_1\setminus S) \times B|s_3^{(n)})-\Psi((\S_1\setminus S) \times B|s_3)\right),
\end{aligned}
\]
where the first equality holds because $\{S,\S_1\setminus S\}$ is a partition of  $\S_1,$ the  inequality follows from the sub-additive properties of upper limits and supremums, and the last  equality holds because the marginal kernel $\Psi(\S_1,\,\cdot\,|\,\cdot\,)$ on $\S_2$ given $\S_3$ is continuous in total variation. So, inequality \eqref{eq:equivWTV0} for arbitrary
open set $\oo\subset\S_1$ follows from inequality \eqref{eq:equivWTV1} for a closed set $C=\S_1\setminus \oo.$ Vice versa, inequality \eqref{eq:equivWTV1} for arbitrary
closed set $\oo\subset\S_1$ follows from inequality \eqref{eq:equivWTV0} for an open set $\oo=\S_1\setminus C.$ That is, (b) $\Leftrightarrow$ ({c}).

({c}) $\Rightarrow$ ({d}). Let $s_3^{(n)}$ converge to $s_3$ in $\S_3,$ and let $A\in\B(\S_1)$ be such that $\Psi(\partial A,\S_2|s_3)=0.$
We shall prove \eqref{eq:equivWTV2}. Indeed, since $\Psi((\bar{A}\setminus A^o)\times \S_2|s_3)=\Psi(\partial A\times\S_2|s_3)=0,$ we have that $\Psi(A^o \times B|s_3)=\Psi(A \times B|s_3)=\Psi(\bar{A} \times B|s_3)$ for each $B\in\B(\S_2).$ Moreover, since $A^o\subset A\subset \bar{A}$ and (b) $\Leftrightarrow$ ({c}), then inequality
\eqref{eq:equivWTV0} applied to $\oo=A^o$ and inequality \eqref{eq:equivWTV1} applied to $C=\bar{A}$ imply
\[
\begin{aligned}
0\le \ilim_{n\to\infty}& \inf_{B\in \B(\S_2)\setminus\{\emptyset\}} \left( \Psi(A^o \times B|s_3^{(n)})-\Psi(A^o \times B|s_3)\right)
\\ &\le \ilim_{n\to\infty} \inf_{B\in \B(\S_2)\setminus\{\emptyset\}} \left( \Psi(A \times B|s_3^{(n)})-\Psi(A \times B|s_3)\right)\\
\le \slim_{n\to\infty}& \sup_{B\in \B(\S_2)\setminus\{\emptyset\}} \left( \Psi(A \times B|s_3^{(n)})-\Psi(A \times B|s_3)\right)
\\ &\le\slim_{n\to\infty} \sup_{B\in \B(\S_2)\setminus\{\emptyset\}} \left( \Psi(\bar{A} \times B|s_3^{(n)})-\Psi(\bar{A} \times B|s_3)\right)\le 0,
\end{aligned}
\]
that is, \eqref{eq:equivWTV2} holds because $\Psi(S \times \emptyset|s)=0$ for each $S\in\B(\S_1)$ and $s\in\S_3.$

({d}) $\Rightarrow$ ({a}). Let (d) hold. Let $s_3^{(n)}$ converge to $s_3$ in $\S_3,$ and let $f$ be a bounded continuous function on $\S_1.$
We shall prove \eqref{eq:equivWTV3}. Indeed, similarly to Parthasarathy~\cite[pp.~41-42]{Part}, let us set
\[
\Psi_f(S,\S_2|s_3):=\Psi(\{s_1\in\S_1\,:\,f(s_1)\in S\},\S_2|s_3),\quad S\in\B(\mathbb{R}).
\]
Since $f$ is a bounded function, there exists a bounded interval $(a,b)$ such that
$a<f(s_1)<b$ for each $s_1\in\S_1,$ and $\Psi_f(\,\cdot\,,\S_2|s_3)$ is concentrated on $(a,b).$  Moreover, the set $\{s\in\mathbb{R}\,:\,\Psi_f(\{s\},\S_2|s_3)>0\}$ is countable or finite. Therefore, for a fixed $\varepsilon >0$ there exist $N_\varepsilon=1,2,\ldots$ and $t_\varepsilon^{(0)}=a<t_\varepsilon^{(1)}<t_\varepsilon^{(2)}<\ldots<t_\varepsilon^{(N_\varepsilon)}=b$ such that $t_\varepsilon^{(i)}-t_\varepsilon^{(i-1)}<\varepsilon$ and $\Psi_f(\{s_1\in\S_1\,:\,f(s_1)=t_\varepsilon^{(i)}\},\S_2|s_3)=0$ for each $i=1,2,\ldots,N_\varepsilon.$

Consider the family of disjoint sets  $\{A^{(i)}:=\{s_1\in\S_1\,:\,t_\varepsilon^{(i-1)}\le f(s_1)<t_\varepsilon^{(i)}\}\}_{i=1}^{N_\varepsilon}.$ Note that $\S_1=\cup_{i=1}^{N_\varepsilon}A^{(i)}.$ Moreover, since $\partial A^{(i)}\subset \{s_1\in\S_1\,:\,f(s_1)=t_\varepsilon^{(i-1)}\}\cup\{s_1\in\S_1\,:\,f(s_1)=t_\varepsilon^{(i)}\},$ we have that $\Psi(\partial A^{(i)},\S_2|s_3)=0,$ and therefore \eqref{eq:equivWTV2} holds with $A=A^{(i)}$ for each $i=1,2,\ldots,N_{\varepsilon}.$ Consequently, for  $f_\varepsilon(s_1):=\sum_{i=1}^{N_\varepsilon}t_{i-1}\h\{s_1\in A^{(i)}\},$ $s_1\in\S_1,$ and for each $n=1,2,\ldots,$
\begin{equation}\label{eqIII}
\begin{aligned}
\sup_{B\in \B(\S_2)}& \left| \int_{\S_1} f(s_1) \Psi(ds_1,B|s_3^{(n)})-\int_{\S_1} f(s_1) \Psi(ds_1,B|s_3)\right|\\
\le&{\color{black} I_1^{(n,\varepsilon)}+ I_2^{(n,\varepsilon)}+ I_3^{(n,\varepsilon)}}\\
\le& 2\varepsilon + \sum_{i=1}^{N_\varepsilon} |t_\varepsilon^{(i-1)}|\sup_{B\in \B(\S_2)} | \Psi(A^{(i)} \times B|s_3^{(n)})-\Psi(A^{(i)} \times B|s_3)|,
\end{aligned}
\end{equation}
where
\[
 {\color{black}
\begin{aligned}
 I_1^{(n,\varepsilon)}&:=\sup_{B\in \B(\S_2)}  \int_{\S_1} \left|f(s_1) -f_\varepsilon(s_1)\right| \Psi(ds_1,B|s_3^{(n)}),\\
 I_2^{(n,\varepsilon)}&:=\sup_{B\in \B(\S_2)}  \int_{\S_1} | f(s_1)-f_\varepsilon(s_1) | \Psi(ds_1,B|s_3)\\
 I_3^{(n,\varepsilon)}&:=\sup_{B\in \B(\S_2)}\left| \int_{\S_1} f_\varepsilon(s_1) \Psi(ds_1,B|s_3^{(n)})-\int_{\S_1} f_\varepsilon(s_1) \Psi(ds_1,B|s_3)\right|,
 \end{aligned}}
\]
{\color{black}and the second inequality in \eqref{eqIII}} holds because $|f(s_1)-f_\varepsilon(s_1)|<\varepsilon$ for each $s_1\in\S_1.$ Letting $n\to\infty,$
\[
\slim_{n\to\infty}\sup_{B\in \B(\S_2)} \left| \int_{\S_1} f(s_1) \Psi(ds_1,B|s_3^{(n)})-\int_{\S_1} f(s_1) \Psi(ds_1,B|s_3)\right|\le 2\varepsilon
\]
because \eqref{eq:equivWTV2} holds with $A=A^{(i)},$ $i=1,2,\ldots,N_\varepsilon.$ Since $\varepsilon>0$ is an arbitrary, \eqref{eq:equivWTV3} holds.$~$\hfill$\Box$

Before the proof of Theorem~\ref{th:concept} we provide an auxiliary lemma.
This lemma is a version of Lemma~5.2 from Feinberg et al. \cite{FKZ} for the class of stochastic kernels satisfying Assumption~\ref{AssKern}.
\begin{lemma}\label{lm:setminuscup}
Let Assumption~\ref{AssKern} hold, and let an arbitrary $s_3\in\S_3$ be fixed. Then
for each $\tilde{\oo}\in \tau_b^{s_3}(\S_1)$ and for each finite union $\oo=\cup_{i=1}^ k {\oo}_{i},$ $k=1,2,\ldots,$ of sets
$\oo_{i}\in\tau_b^{s_3}(\S_1),$ $i=1,2,\ldots,k,$ the set of functions $\fff^\Psi_{\tilde{\oo}\setminus \oo}$ is equicontinuous at $s_3.$
\end{lemma}
\underline{Proof.}
Let $\A^{k}: = \{\cap_{m = 1}^{k} \oo_{j_m}: 1\le j_1\le j_2\le \ldots\le j_{k}\le k\},$ $k=1,2,\ldots,$ be the finite set of all possible intersections of the elements of the tuple  $\{\oo_1,\oo_2, \ldots, \oo_{k}\},$ and let $\hat{\A}^k:=\A^k \cup \{\S_1\}$ be the finite set obtained by adding the single element $\S_1$ to $\A^k.$ Assumption~\ref{AssKern}(ii) imply that the sets of functions $\fff^\Psi_{\tilde{\oo}\setminus \oo}$ is equicontinuous at $s_3$ because
\[
\begin{aligned}
\sup\limits_{B\in \B(\S_2)}&|\Psi ((\tilde{\oo}\setminus \oo)\times B  | s_3') - \Psi ((\tilde{\oo}\setminus \oo)\times B  | s_3)| \\
\le&\sup\limits_{B\in \B(\S_2)} |\Psi (\tilde{\oo}\times B  | s_3') - \Psi (\tilde{\oo}\times B  | s_3)|\end{aligned}
\]
\[
\begin{aligned}
&+\sup\limits_{B\in \B(\S_2)} |\Psi ((\tilde{\oo}\cap \oo)\times B  | s_3') - \Psi ((\tilde{\oo}\cap \oo)\times B  | s_3)|\\
\le& \sum_{D \in \hat{\A}^k} \sup\limits_{B\in \B(\S_2)}|\Psi ((\tilde{\oo}\cap D)\times B | s_3') - \Psi ((\tilde{\oo}\cap D)\times B| s_3)|
\to 0,
\end{aligned}
\]
as $s_3'\to s_3,$ where the first inequality holds because $\tilde{\oo}= (\tilde{\oo}\setminus \cal{O})\cup (\tilde{\oo}\cap \oo)$ and $(\tilde{\cal{O}}\setminus \cal{O})\cap (\tilde{\oo}\cap \oo)=\emptyset,$ and the second inequality follows from the principle of inclusion-exclusion applied to the set $\cal{O}.$$~$\hfill$\Box$

\underline{Proof of Theorem~\ref{th:concept}.}
In view of Theorem~\ref{th:equivWTV}(a,b), it is sufficient to prove that Assumption~\ref{AssKern} holds if and only if the stochastic kernel $\Psi$ on $\S_1\times\S_2$ given $\S_3$ is WTV-continuous.

\textit{Necessity.} Fix an arbitrary $s_3\in\S_3.$ For the topology on $\S_1,$ let us construct its  countable base $\tau_b^{s_3}(\S_1)$
satisfying conditions~(i) and (ii) from Assumption~\ref{AssKern}. For this purpose we firstly note that
every open ball ${\color{black}B(o;\delta)},$ where $\delta>0$ and $o\in\S_1,$ contains open balls
$B{\color{black}(o;{\triangle_o^\delta(i)})}$, $0<\triangle_o^\delta(i)\le\delta$, $i=1,2,\ldots,$ such that
\begin{equation}\label{eq:wtvsc2a}
\triangle_o^\delta(i)\uparrow\delta\quad\mbox{as}\quad i\to\infty,
\end{equation}
and
\begin{equation}\label{eq:wtvsc2}
{\color{black}\Psi((\bar{B}(o;{\triangle_o^\delta(i)})\setminus B(o;{\triangle_o^\delta(i)}))\times\S_2|s_3)=\Psi(S(o;{\triangle_o^\delta(i)})\times\S_2|s_3)=0,}
\end{equation}
that is, {\color{black}$B(o;{\triangle_o^\delta(i)})$} is a continuity set for the probability measure $\Psi(\,\cdot\,|s_3)$ for each $i=1,2,\ldots;$
Parthasarathy \cite[p.~50]{Part}.

Secondly, we set {\color{black}$\oo:=\cap_{j=1}^kB(o_j;{\triangle_{o_j}^{\delta_j}(i_j)})$ and $\tilde{\oo}:= \cap_{j=1}^k\bar{B}(o_j;{\triangle_{o_j}^{\delta_j}(i_j)})$}  for a some natural number $k=1,2,\ldots,$ for  a finite sequence of natural numbers $i_1,i_2,\ldots,i_k$ for a finite sequence of points $o_1,o_2,\ldots,o_k$ from $\S_1,$ and for a finite sequence of positive constants $\delta_1,\delta_2,\ldots,\delta_k.$ We observe that $\partial\oo=\tilde{\oo}\setminus\oo.$ Let us prove that the set of functions $\fff^\Psi_{\oo},$ defined in \eqref{eq:familyoffunctions} with $A=\oo,$ is equicontinuous at~$s_3.$ Indeed, since the stochastic kernel $\Psi$ on $\S_1\times\S_2$ given $\S_3$ is uniform semi-Feller,
equality \eqref{eq:equivWTV0new1} follows from Theorem~\ref{th:equivWTV}(a,d) because
\begin{equation*}
0\le \Psi((\tilde{\oo}\setminus\oo)\times B|s_3)\le \Psi((\tilde{\oo}\setminus\oo)\times\S_2|s_3)=0
\end{equation*}
for each $B\in\B(\S_2),$ where the second inequality holds because $(\tilde{\oo}\setminus\oo)\times B\subset (\tilde{\oo}\setminus\oo)\times\S_2$ for each $B\in\B(\S_2),$ and
the equality holds because $(\tilde{\oo}\setminus\oo)\times\S_2\subset (\cup_{j=1}^kS{\color{black}(o_j; {\triangle_{o_j}^{\delta_j}(i_j)})})\times\S_2$ and $\Psi(S{\color{black}(o_j;{\triangle_{o_j}^{\delta_j}(i_j)})}\times\S_2|s_3)=0$ for all $j=1,2,\ldots,k.$

Finally, according to Rudin \cite[Exercise~2.11]{Rudin}, since the metric space $\S_1$ is separable,
there exists a sequence $\{s^{(j)}\}_{j=1,2,\ldots}\subset \oo$ such
that the set $\{{\color{black}B(s^{(j)};\delta)}\,:\, \delta \in\mathbb{Q}_{>0},\, j=1,2,\ldots\}$ is a countable base of the topology on $\S_1,$ where $\mathbb{Q}_{>0}$ is the set of positive rational numbers. {\color{black} The} set
{\color{black}$\tau_b^{s_3}(\S_1):=\{B({\color{black}s^{(j)};\triangle_{s^{(j)}}^{\delta(i)}})\,:\, \delta \in\mathbb{Q}_{>0},\, i,j=1,2,\ldots\}\cup\{\S_1\}$} is a countable base of the topology on $\S_1$ because, according to \eqref{eq:wtvsc2a}, {\color{black} $B({\color{black}s^{(j)};\delta})=\mathop{\cup}\limits_{i=1,2,\ldots}B({\color{black}{s^{(j)};\triangle_{s^{(j)}}^\delta(i)}})$} for each $j=1,2,\ldots$ and $\delta\in\mathbb{Q}_{>0}.$ Moreover, for each finite intersection $\hat{\oo}=\cap_{i=1}^ k {\oo}_{i}$ of sets
$\oo_{i}\in\tau_b^{s_3}(\S_1),$ $i=1,2,\ldots,k,$
the set of
functions $\fff^\Psi_{\hat{\oo}}$ is equicontinuous at $s_3,$ where set of
functions $\fff^\Psi_{\S_1}$ is equicontinuous at $s_3$ because the marginal kernel $\Psi(\S_1,\,\cdot\,|\,\cdot\,)$ on $\S_2$  given $\S_3$ is continuous in total variation, that is, Assumption~\ref{AssKern} holds.

\textit{Sufficiency.} Assumption~\ref{AssKern} implies that the marginal kernel $\Psi(\S_1,\,\cdot\,|\,\cdot\,)$ on $\S_2$  given $\S_3$ is continuous in total variation because, by the definition, equicontinuity of the set $\fff^\Psi_{\S_1}$ at a point $s_3\in\S_3$ is equivalent to the continuity in total variation of the marginal kernel $\Psi(\S_1,\,\cdot\,|\,\cdot\,)$ on $\S_2$ given $\S_3$ at the point $s_3.$

Let us prove the WTV-continuity of the stochastic kernel $\Psi$ on $\S_1\times\S_2$ given $\S_3.$ For this purpose we fix an arbitrary element $s_3\in\S_3$ and a sequence $\{s_3^{(n)}\}_{n=1,2,\ldots}\subset \S_3$ such that $s_3^{(n)}\to s_3$ as $n\to\infty.$ Let us prove that (\ref{eq:wtvsc}) holds for an arbitrary fixed $\oo\in\tau(\S_1).$ Indeed, Assumption~\ref{AssKern}(ii) implies the existence of a tuple $\{\oo_1,\oo_2,\ldots\}\subset \tau_b^{s_3}(\S_1)$ such that $\oo=\cup_{j=1}^\infty\oo_j.$ Setting
$A_k:=\cup_{j=1}^k\oo_j,$ $k=1,2,\ldots,$ and $A_0:=\emptyset,$ we note that Lemma~\ref{lm:setminuscup} implies that the set of
functions $\fff^\Psi_{\oo_k\setminus A_{k-1}}$ is equicontinuous at $s_3$ for each $k=1,2,\ldots\,.$
Thus,
\begin{align*}
&\ilim_{n\to\infty} \inf_{B\in \B(\S_2)} \left( \Psi(\oo \times B|s_3^{(n)})-\Psi(\oo \times B|s_3)\right)\\
&=\ilim_{n\to\infty} \inf_{B\in \B(\S_2)} \left( \Psi((\cup_{j=1}^\infty (A_k\setminus A_{k-1})) \times B|s_3^{(n)})\right.\\ &\left.-\Psi((\cup_{j=1}^\infty (A_k\setminus A_{k-1})) \times B|s_3)\right)\\
&=\ilim_{n\to\infty} \inf_{B\in \B(\S_2)} \sum_{k=1}^\infty\left( \Psi((A_k\setminus A_{k-1}) \times B|s_3^{(n)})-\Psi((A_k\setminus A_{k-1}) \times B|s_3)\right)
\\
&\ge \ilim_{n\to\infty} \sum_{k=1}^\infty \inf_{B\in \B(\S_2)} \left( \Psi((A_k\setminus A_{k-1}) \times B|s_3^{(n)})-\Psi((A_k\setminus A_{k-1}) \times B|s_3)\right)\\
&\ge \sum_{k=1}^\infty \ilim_{n\to\infty} \inf_{B\in \B(\S_2)} \left( \Psi((A_k\setminus A_{k-1}) \times B|s_3^{(n)})-\Psi((A_k\setminus A_{k-1}) \times B|s_3)\right)= 0,
\end{align*}
where first two equalities hold because $\oo=\cup_{j=1}^\infty\oo_j=\cup_{j=1}^\infty (A_k\setminus A_{k-1})$ and $(A_j\setminus A_{j-1})\cap (A_i\setminus A_{i-1})=\emptyset$ for each $i\ne j,$ the first inequality follows from the basic property of infimums, the second inequality follows from Fatou's lemma because each summand is bounded below by $-1$ since for each $k=1,2,\ldots$
\[
\begin{aligned}
\inf_{n=1,2,\ldots}\inf_{B\in \B(\S_2)}& \left( \Psi((A_k\setminus A_{k-1}) \times B|s_3^{(n)})-\Psi((A_k\setminus A_{k-1}) \times B|s_3)\right)\\ &\ge -\Psi((A_k\setminus A_{k-1}) \times \S_2|s_3),
\end{aligned}
\]
and $\sum_{k=1}^\infty \Psi((A_k\setminus A_{k-1}) \times \S_2|s_3)=1,$ and the last equality holds because the set of
functions $\fff^\Psi_{A_k\setminus A_{k-1}}=\fff^\Psi_{\oo_k\setminus A_{k-1}}$ is equicontinuous at $s_3$ for each $k=1,2,\ldots\,.$ Therefore, inequality (\ref{eq:wtvsc}) holds for each $\oo\in\tau(\S_1),$ that is, the stochastic kernel $\Psi$ on $\S_1\times\S_2$ given $\S_3$ is WTV-continuous.$~$\hfill$\Box$

\underline{Proof of Theorem~\ref{th:extra}.}
\textit{Sufficiency.} Let  $(s_3^{(n)},s_4^{(n)})\to (s_3,s_4)$ in $\S_3\times\S_4$ as $n\to\infty.$ Consider the sequence of
probability measures $\{\mu^{(n)},\mu\}_{n=1,2,\ldots}$ such that $\mu^{(n)}(C)=\h\{s_3^{(n)}\in C\}$ and $\mu(C)=\h\{s_3\in C\}$ for each
$C\in\B(\S_3)$ and $n=1,2,\ldots\,.$ Since $(\mu^{(n)})_{n=1,2,\ldots}$ converges weakly to $\mu,$ and the stochastic kernel
$\Psii$ on $\S_1\times\S_2$ given $\P(\S_3)\times\S_4$ is WTV-continuous, we obtain that
\[
\lim_{n\to\infty} \inf_{B\in \B(\S_2)} \left( \Xi(\oo \times B|s_3^{(n)},s_4^{(n)})-\Xi(\oo \times B|s_3,s_4)\right)= 0
\]
for each $\oo \in\tau (\S_1),$ that is, the stochastic kernel $\Xi$ on $\S_1\times\S_2$ given $\S_3\times\S_4$ is WTV-continuous.

\textit{Necessity.} {\color{black} Semi-uniform Fellerness} of the stochastic kernel $\Xi$ on $\S_1\times\S_2$ given $\S_3\times\S_4$ implies that for each $\oo \in\tau (\S_1)$ the
set of functions
$\fff_\oo^\Xi=\{  (s_4,s_3)\mapsto \Xi(\oo\times B |s_3,s_4):\, B\in \B(\S_2)\}$
is lower semi-equicontinuous. Theorem~\ref{th:PrLESC} applied to $S_1:=\S_4,$ $S_2:=\S_3,$ $S_3:=\P(\S_3),$ $\mathcal{A}:=\fff_\oo^\Xi,$ and $\psi(\,\cdot\,|\mu):=\mu(\,\cdot\,)$ for $\mu\in S_3,$
implies that the set of functions $\{(\mu,s_4)\mapsto \Psii(\oo\times B|\mu,s_4)\,:\, B\in\B(\S_2)\}$ is lower semi-equicontinuous because the stochastic kernel $\psi$ on $S_2$ given $S_3$ is weakly continuous. Since $\oo \in\tau (\S_1)$ is an arbitrary, the stochastic kernel $\Psii$ on
$\S_1\times\S_2$ given $\P(\S_3)\times\S_4$ is {\color{black}semi-uniform Feller}.$~$\hfill$\Box$

\section*{Acknowledgements} We thank Janey (Huizhen) Yu {\color{black} and Yi Zhang} for valuable remarks. The second and the third authors were partially supported by the National Research Foundation of Ukraine, Grant No.~2020.01/0283.  The second author was partially supported by a U4U non-residential fellowship by UC Berkeley Economics/Haas.


\begin{thebibliography}{}
%
%


\bibitem{Ao} Aoki, M. (1965) Optimal control of partially observable Markovian systems. \emph{J. Franklin Inst.} 280(5): 367--386.

\bibitem{As}
 {\AA}str\"om, K.J. (1965). Optimal control of Markov processes with incomplete state information. \emph{Journal of Mathematical Analysis and Applications} 10: 174--205.

{ \color{black}
 \bibitem{BR} B\"auerle, N., Rieder, U. (2011) \textit{Markov Decision Processes with Applications to Finance,} Springer-Verlag, Berlin.
 }

\bibitem{BS} Bertsekas, D.P.,  Shreve S.E. (1978) \textit{Stochastic Optimal Control: The Discrete-Time
Case,} Academic Press, New York 

\bibitem{bogachev} Bogachev, V.I. (2007) \textit{Measure Theory, Volume II,}
Springer-Verlag, Berlin.


\bibitem{Dy}  Dynkin, E.B. (1965) Controlled random sequences. \emph{Theory
Probab. Appl.} 10(1): 1--14.

{\color{black}
\bibitem{DY} Dynkin, E.B.,   Yushkevich A.A. (1979) \textit{{C}ontrolled
{M}arkov {P}rocesses,} Springer-Verlag, New York.}

{\color{black}
\bibitem{FKZ22} Feinberg, E.A.,  Kasyanov, P.O. (2022) Equivalent conditions for weak continuity of nonlinear filters, arXiv2207.07544.
}

\bibitem{FKL2} Feinberg, E.A.,  Kasyanov, P.O., Liang, Y. (2020) Fatou's lemma in its classical form and Lebesgue's
convergence theorems for varying measures with applications to Markov decision processes. \emph{Theory Probab. Appl.} 65(2): 270--291.

{\color{black}\bibitem{FKR} Feinberg, E.A.,  Kasyanov, P.O., Royset, J.O. (2023)   Epi-convergence of expectation functions under varying measures and integrands. \emph{J. Convex Anal.} 30(2) (to appear) arXiv:2208.03805}

{\color{black}
\bibitem{FKNMOR}  Feinberg, E.A.,  Kasyanov, P.O., Zadoianchuk, N.V. (2012)   Average-cost Markov decision
processes with weakly continuous transition probabilities.
\emph{Math. Oper. Res.} 37(4): 591--607.
}

\bibitem{Steklov} Feinberg, E.A., Kasyanov, P.O., Zgurovsky, M.Z. (2014) Convergence of probability measures and Markov decision models with incomplete information,
\emph{Proceedings of the Steklov Institute of Mathematics,} 287 (1), 96--117.

\bibitem{UFL} Feinberg, E.A.,  Kasyanov, P.O., Zgurovsky, M.Z. (2016) Uniform Fatou's lemma,  \emph{Journal of Mathematical Analysis and Applications}, 444(1), 550--567.

\bibitem{FKZ} Feinberg, E.A., Kasyanov, P.O., Zgurovsky, M.Z. (2016) Partially observable total-cost Markov decision processes with weakly continuous transition probabilities, \emph{Math. Oper. Res.}, 41(2), 656--681.

\bibitem{FKZSIAM} Feinberg, E.A.,  Kasyanov, P.O., Zgurovsky, M.Z. (2022) Markov decision processes with incomplete information and semi-uniform Feller transition probabilities,  \emph{SIAM Journal on Control and Optimization} 60: 2488--2513.

\bibitem{FKZConf} Feinberg, E.A.,  Kasyanov, P.O., Zgurovsky, M.Z. (2021) Average cost Markov decision processes with  semi-uniform Feller transition probabilities. Piunovskiy, A., Zhang, Yi, eds.  \emph{ Modern Trends in Controlled Stochastic Processes:
Theory and Applications,} Springer Nature, Cham. 



{\color{black}
\bibitem{HL}  Hern\'{a}ndez-Lerma, O. (1989) \textit{Adaptive Markov Control Processes,} Springer-Verlag, New York.
}

\bibitem{Saldi} Kara, A.D., Saldi, N., Y\"uksel, S. (2019) Weak Feller property of non-linear filters. \emph{Systems $\&$ Control Letters} 134: 104512.

\bibitem{Ma} Ma, L. (2021) Sequential convergence on the space of Borel measures. arXiv:2102.05840

\bibitem{Papa} Papanicolaou, G.C. (1978) Asymptotic analysis of stochastic equations. Rosenblatt M., ed.  \emph{Studies in Probability Theory}  Mathematical Association of America, Washington DC, 111--179.

\bibitem{Part}  Parthasarathy, K.R. (1967)  \textit{Probability Measures on Metric Spaces,} Academic Press, New York.

{\color{black}
\bibitem{Plat} Platzman, L.K. (1980) Optimal infinite-horizon undiscounted control of finite probabilistic systems. \emph{SIAM Journal on Control and Optimization} 18(4): 362--380.

\bibitem{Rh}  Rhenius, D. (1974) Incomplete information in Markovian decision models. \emph{Ann. Statist.} 2(6): 1327--1334.

{\color{black}\bibitem{RW} Rockafellar, R. T., Wets, R. J.-B. (1998) \textit{Variational Analysis}, Springer, Berlin.}

\bibitem{Ri} Rieder, U. (1975) Bayesian dynamic programming. \emph{Adv. Appl.
Probab.}  7(2): 330--348.
}



\bibitem{Rudin} Rudin, W. (1964) \textit{Principles of Mathematical Analysis,} Second edition, McGraw-Hill Inc, New York.

{\color{black}
\bibitem{RS} Runggaldier, W.J.,   Stettner, L. (1994)
\emph{Approximations of Discrete Time Partially Observed Control Problems,}
Applied Mathematics Monographs CNR, Giardini Editori, Pisa.}
\bibitem{SLY} Saldi, N., T. Linder, T., Yuksel, (2018)  \emph{Finite Approximations in Discrete-Time Stochastic
Control: Quantized Models and Asymptotic Optimality,} Springer, Cham.


{\color{black}\bibitem{sha} Sch\"al, M. (1975) On dynamic programming: compactness of the space of policies. \emph{Stoch. Process.  Appl.} 3:345--364.}

\bibitem{Shi64}{\color{black}
Shiryaev, A.N. (1964) On the theory of decision functions and control by an observation process with incomplete data. \textit{Transactions of the Third Prague
Conference on Information Theory, Statistical Decision Functions,
Random Processes} (Prague, 1962), pp.~657-681 (in Russian); Engl.
transl. in \textit{Select. Transl. Math. Statist. Probab.}
6(1966), 162-188.}


\bibitem{Shi}
Shiryaev, A.N. (1967) Some new results in the theory of controlled
random processes. \textit{Transactions of the Fourth Prague
Conference on Information Theory, Statistical Decision Functions,
Random Processes} (Prague, 1965), pp.~131-201 (in Russian); Engl.
transl. in \textit{Select. Transl. Math. Statist. Probab.}
8(1969), 49-130.

\bibitem{SS} Smallwood, R.D., Sondik E.J. (1973) The optimal control of partially
observable Markov processes over a finite horizon. \emph{Oper.
Res.} 21(5): 1071--1088.

\bibitem{Yu} Yushkevich AA (1976) Reduction of a controlled Markov model with incomplete data to a problem with complete information
in the case of Borel state and control spaces. \emph{Theory
Probab. Appl.} 21(1): 153--158.
\end{thebibliography}
\end{document}